\documentclass[numreferences]{kluwer}    

\input amssym.def

  \newtheorem{thm}{THEOREM}
  \newtheorem{cor}[thm]{COROLLARY}

  \newtheorem{REM}[thm]{REMARK}

  \rmtheorem{alg}
  \newtheorem{defn}[thm]{DEFINITION}
  \rmtheorem{defn}
  \newtheorem{exmp}[thm]{EXAMPLE}
  \rmtheorem{exmp}

\def\complessi{{\mathchoice {\setbox0=\hbox{$\displaystyle\rm C$}\hbox{\hbox
to0pt{\kern0.4\wd0\vrule height0.9\ht0\hss}\box0}}
{\setbox0=\hbox{$\textstyle\rm C$}\hbox{\hbox
to0pt{\kern0.4\wd0\vrule height0.9\ht0\hss}\box0}}
{\setbox0=\hbox{$\scriptstyle\rm C$}\hbox{\hbox
to0pt{\kern0.4\wd0\vrule height0.9\ht0\hss}\box0}}
{\setbox0=\hbox{$\scriptscriptstyle\rm C$}\hbox{\hbox
to0pt{\kern0.4\wd0\vrule height0.9\ht0\hss}\box0}}}}
\def\bbbe{{\mathchoice {\setbox0=\hbox{\smalletextfont e}\hbox{\raise
0.1\ht0\hbox to0pt{\kern0.4\wd0\vrule width0.3pt height0.7\ht0\hss}\box0}}
{\setbox0=\hbox{\smalletextfont e}\hbox{\raise
0.1\ht0\hbox to0pt{\kern0.4\wd0\vrule width0.3pt height0.7\ht0\hss}\box0}}
{\setbox0=\hbox{\smallescriptfont e}\hbox{\raise
0.1\ht0\hbox to0pt{\kern0.5\wd0\vrule width0.2pt height0.7\ht0\hss}\box0}}
{\setbox0=\hbox{\smallescriptscriptfont e}\hbox{\raise
0.1\ht0\hbox to0pt{\kern0.4\wd0\vrule width0.2pt height0.7\ht0\hss}\box0}}}}
\def\razionali{{\mathchoice {\setbox0=\hbox{$\displaystyle\rm Q$}\hbox{\raise
0.15\ht0\hbox to0pt{\kern0.4\wd0\vrule height0.8\ht0\hss}\box0}}
{\setbox0=\hbox{$\textstyle\rm Q$}\hbox{\raise
0.15\ht0\hbox to0pt{\kern0.4\wd0\vrule height0.8\ht0\hss}\box0}}
{\setbox0=\hbox{$\scriptstyle\rm Q$}\hbox{\raise
0.15\ht0\hbox to0pt{\kern0.4\wd0\vrule height0.7\ht0\hss}\box0}}
{\setbox0=\hbox{$\scriptscriptstyle\rm Q$}\hbox{\raise
0.15\ht0\hbox to0pt{\kern0.4\wd0\vrule height0.7\ht0\hss}\box0}}}}
\def\bbbt{{\mathchoice {\setbox0=\hbox{$\displaystyle\rm
T$}\hbox{\hbox to0pt{\kern0.3\wd0\vrule height0.9\ht0\hss}\box0}}
{\setbox0=\hbox{$\textstyle\rm T$}\hbox{\hbox
to0pt{\kern0.3\wd0\vrule height0.9\ht0\hss}\box0}}
{\setbox0=\hbox{$\scriptstyle\rm T$}\hbox{\hbox
to0pt{\kern0.3\wd0\vrule height0.9\ht0\hss}\box0}}
{\setbox0=\hbox{$\scriptscriptstyle\rm T$}\hbox{\hbox
to0pt{\kern0.3\wd0\vrule height0.9\ht0\hss}\box0}}}}
\def\bbbs{{\mathchoice
{\setbox0=\hbox{$\displaystyle \rm S$}\hbox{\raise0.5\ht0\hbox
to0pt{\kern0.35\wd0\vrule height0.45\ht0\hss}\hbox
to0pt{\kern0.55\wd0\vrule height0.5\ht0\hss}\box0}}
{\setbox0=\hbox{$\textstyle \rm S$}\hbox{\raise0.5\ht0\hbox
to0pt{\kern0.35\wd0\vrule height0.45\ht0\hss}\hbox
to0pt{\kern0.55\wd0\vrule height0.5\ht0\hss}\box0}}
{\setbox0=\hbox{$\scriptstyle \rm S$}\hbox{\raise0.5\ht0\hbox
to0pt{\kern0.35\wd0\vrule height0.45\ht0\hss}\raise0.05\ht0\hbox
to0pt{\kern0.5\wd0\vrule height0.45\ht0\hss}\box0}}
{\setbox0=\hbox{$\scriptscriptstyle\rm S$}\hbox{\raise0.5\ht0\hbox
to0pt{\kern0.4\wd0\vrule height0.45\ht0\hss}\raise0.05\ht0\hbox
to0pt{\kern0.55\wd0\vrule height0.45\ht0\hss}\box0}}}}

\def\d{\,{\rm d}} 

\def\supp{{\rm supp}}
\def\phi{{\varphi}}
\def\epsilon{{\varepsilon}}
\def\aut{{\rm AUT}}

\def\reali{{\Bbb R}}
\def\complessi{{\Bbb C}}
\def\zzz{{\Bbb Z}}
\def\naturali{{\Bbb N}}
\def\boxf{{$\ulcorner \! \! \lrcorner \! \! \! \! \llcorner \! \! \urcorner$}}
\def\QED{{\ifmmode\sq\else{\unskip\nobreak\hfil
\penalty50\hskip1em\null\nobreak\hfil{\boxf}
\parfillskip=0pt\finalhyphendemerits=0\endgraf}\fi} \vskip 12 pt} 

\def\parti#1{{{\cal P} \left({#1}\right )}}
\def\har{{\cal H}}

\def\calM{{\cal M}}

\begin{document}                                                                                   
\begin{article}
\begin{opening}         

\title{Strong and weak mean value properties on trees}

\author{Fabio Zucca}
\institute{Universit\'e de Cergy-Pontoise \\ 
	D\'epartement de Math\'ematiques \\
	2 rue Adolphe Chauvin, 95302 Pontoise, France.}

\runningtitle{Strong and weak mean value properties on trees}
\runningauthor{Fabio Zucca}

\begin{abstract}
We consider the mean value properties for finite variation measures
with respect to a Markov operator in a discrete environnement. We prove
equivalent conditions for the weak mean value property in the case of general 
Markov operators and for the strong mean value property in the case of 
transient Markov operators adapted to a tree structure. In this
last case, conditions for the equivalence between weak and strong mean
value properties are given. 
\end{abstract}

\classification{AMS 2000 MSC}{60J10, 31C05}

\end{opening}


\begin{section}[1]{Introduction} \label{1}

This paper deals with the discrete analogous of the theory of harmonic functions. The set of (discrete) harmonic functions, 
defined on an at most countable set $X$, can be constructed 
starting from
a Markov operator $P$. 
More precisely, given a Markov operator $P$ acting on real functions defined on
$X$ according to 
$$
(Pf)(x)=\sum_{y\in X} p(x,y)f(y)
$$
(where $p:X \times X \rightarrow [0,1]$ and $\sum_{y\in X} p(x,y)=1$, for
any $x\in X$) provided that $\sum_{y\in X} p(x,y)|f(y)|$ for all $x\in X$,
then we call $f$ a {\it $P$-harmonic function} (or simply a {\it harmonic function}) if $Pf=f$ and we denote the set of harmonic functions by
${\cal H}(X,P)$.

A fundamental tool of this theory is represented by the {\it Martin compactification} 
(see \cite{Woess2}) which allows us to state and solve (under very
general condition) the (discrete) Dirichlet problem. More generally,
using this approach, one can obtain a useful representation of
bounded harmonic functions (see \cite{Cartier} and \cite{Pica-Taible})
${\cal H}^\infty(X,P)$. If $X$ is
a vertex set of a tree and $P$ is {\it adapted} to the graph structure
then a generalization of the previous arguments leads to a representation
of any harmonic function (see \cite{Cartier}).

In this paper we study a mean value property (see Definition~\ref{mean}), which is actually the dual
property of the one commonly known about harmonic functions.

This problem was stated and partially studied in \cite{Zucca1}. The main
results of that paper are on radial trees (with respect to the centre
of symmetry); here we use a completely different technique to 
obtain a generalization 
to the non-symmetric case. In particular
we are interested to the relation between weak and strong MVPs.

Besides we consider the mean walue property with respect to a set of (harmonic) functions
${\cal F}=\{k(\cdot,\xi): \xi \in F\}$ where $F$ is a suitable subset of the martin boundary of 
harmonic measure $1$ (see the next section for the definitions). We show that
there are equivalent conditions to the weak and strong MVP in
terms of MVPs with respect to such a family ${\cal F}$.

We give a brief outline of the paper. In Section~\ref{2} we give some basic definitions and 
we fix the notation. Section~\ref{3} is a short introduction to the main representation 
theorem for harmonic functions on trees (which generalizes, in the case of
trees, the well known Poisson-Martin representation theorem for bounded harmonic functions
(see \cite{Woess2}, Theorem 24.7).
Section~\ref{5} is devoted
to the study of the relation between the weak MVP and the MVP with respect to a suitable 
subset of minimal harmonic functions $k(\cdot,\xi)$ for finite supported signed measure on a general graph
(Theorem~\ref{equiv2}).
The analogous case for the strong MVP on transient trees is considered in Section~\ref{4}
(Theorem~\ref{equiv1}): in particular Corollary~\ref{trees1} guarantees, under
very general conditions, the equivalence between weak and strong MVP on a wide class
of random walk on trees. This propery, which does not hold in general 
(Examples~\ref{ex1} and \ref{ex2}), is useful, since equivalent
conditions for the weak mean value property for bounded variation signed measures 
(with bounded or unbounded support) on general 
graphs are known (see \cite{Zucca1}, Section 6). 
The results of Section~\ref{5} are extended to measures with unbounded support in 
Section~\ref{6}.
Finally, in Section~\ref{7} we briefly discuss some results in the
case of general natural compactifications.

\end{section}

\begin{section}[2]{Preliminaries and basic definitions}
\label{2}


In this section we state some basic definitions which will be 
widely used in the following. 

Given any Markov operator $P$ on $X$ and any probability measure $\lambda_0$ on 
$\parti{X}$ it is well known that, according to Kolmogorov Theorem,
there exists a probability space $(\Omega, {\Sigma}_\Omega, \Bbb P)$ and a (stationary) 
Markov chain $\{Z_n\}_{n\in \naturali}$ (called {\it canonical Markov Chain})
such that $\Bbb P_{Z_0} \equiv \lambda_0$ and 
$\{p(x,y)\}_{x,y \in X}$ represent the transition probabilities of $\{Z_n\}_{n \in \naturali}$.
In particular by $\{Z_n^x\}_{n\in \naturali}$ (given any fixed $x\in X$) we mean the 
canonical Markov chain corresponding to the Dirac initial distribution $\delta_x$.

We denote by $G$ be the generating functions of the transition probabilities
and by $F$ the generating function of the hitting probabilities (see \cite{Woess2}, Chapter 1,
Paragraph 1.B).

When we deal with the graph associated to $(X,P)$ we mean the graph generated by 
the random walk $(X,P)$ (i.~e.~$(X,E_P(X))$ where $E_P:=\{(x,y)\in X\times X: p(x,y) >0\}$).
We consider only {\it irreducible} Markov operators (i.~e.~giving rise to irreducible random walks 
according to the usual definition); hence $(X,E_P)$ is a connected graph. In this
case the properties of recurrence and transience (\cite{Woess2}, Definition 1.14) does not depend on the particular state 
$x\in X$.

Given a topological Hausdorff space $(X, \tau)$, we say that
$(\widehat X, \widehat \tau, i)$ is a compactification of $X$ if and only if
$(\widehat X, \widehat \tau)$ is a compact hausdorff space, $i: X \rightarrow \widehat X$ is
a homeomorphism between $X$ and $i(X)$ with the induced topology and $i(X)$ is dense in
$\widehat X$.
In the case of a countable graph (with
the discrete topology) then $i(X)$ is open and discrete. The boundary of the compactification
is $\partial X:= \widehat X \setminus i(X)$; from now on, we identify $X$ with $i(X)$, hence
the map $i$ will be the inclusion map and the compactification will be denoted simply
by $(\widehat X, \widehat \tau)$ (or shortly by $\widehat X$).

If $\{Z_n\}$ is an irreducible Markov chain
on $X$ and $\widehat X$ is a compactification of 
the state space $X$, then we say that $\{Z_n\}$ converges a.c.~to the boundary
if and only if 
${\rm Pr}(\lim_n Z_n \in \partial X)=1$;
in this case it is well known that $Z_\infty:=\lim_{n \rightarrow +\infty} Z_n$ defines
a.e.~a random variable with values in $\partial X$ (measurable with respect to the Borel
$\sigma$-algebra on $\partial X$). 

\begin{defn}\label{naturalco}
Let $(X,P)$ be a random walk; a compactification $\widehat X$ of 
the discrete set $X$ is called
a {\it natural compactification with respect to $P$} if and only if 
for any $x\in X$ the canonical Markov
process $\{Z_n^x\}$ with initial distribution $\delta_x$ converges 
${\rm Pr}_x$-a.c.~to the boundary $\partial X$.
\end{defn}

It is known (see \cite{Woess2}) that if $\widehat X$ is a natural
compactification, then the $\partial X$-valued random variables 
representing the limits to the boundary of 
$\{Z_n^x\}$ does not depend on $x$; hence there exists
a $\partial X$-valued random variable $Z_\infty$ such that,
for any $x\in X$, we have 
$$
{\rm Pr}_x(\lim_{n \rightarrow +\infty} Z_n^x = Z_\infty)=1.
$$
We define the {\it harmonic measures}
$\nu_x(A):=\Pr_x(Z_\infty^{-1}(A))$ for any borel set $A \subseteq \partial X$ and 
for every $x \in X$. Since
$$
\nu_x=\sum_{y\in X} p(x,y)\nu_y,
$$
it is obvious that these measures are mutually absolutely continuous, whence their
supports coincides and 
$$
\begin{array}{rl}
L^\infty(\partial X,\nu_x) = L^\infty(\partial X, \nu_y) \subseteq 
L^p(\partial X, \nu_y)&=L^p(\partial X, \nu_x), \\
& \forall x,y \in X, \forall p\in [1,+\infty),
\end{array}
$$
where on $\partial X$ we consider the $\sigma$-algebra of Borel.

The {\it Martin
compactification} (see \cite{Woess2} Chapter 4 
Paragraph 24) is constructed by means of the
Martin Kernel $k_o$ 
which is defined by $k_o(x,y):= F(x,y)/F(o,y)$; in particular
it is the unique (up to homeomorphisms) 
smallest compactification of the discrete set $X$ such that 
all the functions $\{k_o(x,\cdot)\}_{x\in X}$ extend continuously to the boundary 
(denoted by ${\cal M}(X,P)$).
The Martin compactification does not depend on the choice of $o$
(in the sense that all the compactifications obtained from different
choices of the reference vertex are homeomorphic);
we denote the extended function again by $k_o(x,\cdot)$ (for any
$x \in X$).
Theorem 24.10 of \cite{Woess2} shows that this is a natural compactification and that,
for any pair of vertices $x,y\in X$, the Radon-Nikodym derivative of
the harmonic measure $\nu_x$ with respect to $\nu_y$ is $k_y(x,\cdot)|_{\partial X}$
(the latter is the extended $\widehat X$-valued Martin kernel restricted
to the boundary $\partial X$).

A sequence of measures $\{\nu_n\}$ on a measurable space $(X,\sigma_\tau)$ 
(where
$\sigma_\tau$ is the Borel $\sigma$-algebra generated by the topology
$\tau$ is said to be {\it weakly convergent} to a measure $\nu$ if
and only if for any $f \in C(X)$,
$$
\int_X f \d \nu_n \ {\buildrel n \rightarrow +\infty \over
\longrightarrow}\ \int_X f \d \nu.
$$
A boundary point $\xi$ of a natural compactification of $(X,P)$
is said to be {\it regular} if and only if, for any sequence 
$\{x_n\}$ in $X$ convergent to $\xi$, we have that the corresponding
sequence of harmonic measures $\{\nu_{x_n}\}$ converges weakly to 
$\delta_\xi$.

\begin{defn}\label{almostregular}
A natural compactification $\widehat X$
of $(X,P)$ is {\it almost regular} if and only if there exists
$x \in X$ ($\Leftrightarrow$ for every $x\in X$) such that
the set of regular point of $\partial X$ has $\nu_x$ measure $1$.
\end{defn}

According to a well known theorem (see Theorem 2.2 of \cite{Woess4} 
for instance) the Dirichlet problem for $P$-harmonic functions is solvable in the 
compactification $\widehat X$ if and only if the compactification 
of $(X,P)$ is natural and every boundary point is regular.

From now on, if not otherwise explicitly stated, the compactification $\widehat X$
of $X$ will be the Martin compactification and the boundary
$\partial X$ will be ${\calM}(X,P)$.

If $X$ is a tree, then for any couple of vertices $x,y \in X$ there 
exists a unique geodesic from $x$ to $y$ denoted by $\Pi[x,y]$. Let us
denote by $X_{x,y}$ the set of vertices visited by the geodesic path
$\Pi[x,y]$; hence for any $x,y,z\in X$,
$X_{x,y}\cap X_{y,z} \cap X_{z,x}$ is a one point set (let us denote it by
$p(x,y,z)$). Since $X_{x,y}=X_{y,x}$ and $X_{x,x}=\{x\}$ then 
$p(x,y,z)$ is invariant under permutation of $x,y,z$ and it is
called the confluent $x \wedge_z y$
of $x$ and $y$ with respect to the reference point
$z$ (see \cite{Woess2} Chapter 1 Paragraph 6.B). Hence $x \wedge_z y=
y \wedge_z x = z \wedge_x y$; moreover the definition
of confluent $\cdot \wedge_o \cdot$
with respect to a reference point $o \in X$ can be extended to the Martin compactification 
$\widehat X$
and for the extended Martin kernel the following relation
$$
k_o(x,\xi)=k_o(x,x \wedge_o \xi), \qquad \forall x \in X, \forall \xi \in \widehat X,
$$
holds
(indeed, in this case, $F(x,y)=F(x, z)F(z,y)$ for any $x,y\in X$ and any $z \in X_{x,y}$).
By continuity and compactness we have that $\xi \mapsto k_o(x, \xi)$ is bounded (and positive)
on $\widehat X$ and
$$
\sup_{y \in X} k_o(x,y) = \max_{\xi \in \widehat X} k_o(x,\xi)=1/F(o,x);
$$
moreover since
$F(x,y)=F(z,y) F(x,x \wedge_z y) / F(z, x \wedge_z y)$
  for any given $x,y,z \in X$,
we have that
$k_o(x,y)=k_o(z,y) k_z(x,y)$ for any $x,y \in \widehat X$ and any $z,o \in X$. 

Given any couple of vertex $o,x \in T$ we call {\it branching subtree} of $x$ with respect to $o$
the set $T_{o,x}:=\{y \in T: x \in \Pi[o,y]\}$ and we denote by $\partial T_{o,x}$
the set $\{\xi \in \partial X: x \in \Pi[o,\xi]\}$.

We turn our attention to the basic definition of mean value property.

\begin{defn}\label{mean}
Let $(X,\Sigma_X)$ be a measurable space,
$o \in X$ and let $\nu$ be a measure on $(X,\Sigma_X)$ with finite variation.
If ${\cal F}$ is a family of functions in
$L^1(|\nu|)$ then we say that $\nu$ has the {\it mean value property} (MVP)  with
respect to ${\cal F}$ and $o$ if
\begin{equation}\label{mean1}
L(h,\nu)(o)=0, \quad \forall h \in {\cal F},
\end{equation}
where
\begin{equation}\label{operatore1}
L(h, \nu)(x):= \int_{X}h \d \nu - \nu(X) h(x), \quad \forall x \in X.
\end{equation}
In particular, if $(X,P)$ is a graph with an adapted random walk (with 
the $\sigma$-algebra of the set af subsets of $X$),
we say that $\nu$ has the {\rm strong mean value property
with respect to $o$}
(resp.~{\rm weak mean value property with respect to $o$}) if
equation~\ref{mean1} holds with ${\cal F}
\equiv \har(X,P)\cap L^1(|\nu|)$  (resp.~${\cal F}
\equiv \har^\infty(X,P):= \har(X,P) \cap
l^\infty(X)$).
If $X$ is a tree with root $o$ then a mean value property will always be
with respect to $o$.
\end{defn}

We stress that 
this one is a property of measures and not 
of functions as in the usual case.

\begin{REM}\label{radial}
If we have an adapted random walk $(X,P)$ on the graph
$X$, $o\in X$ and 
$P$ is $\Gamma_o$-invariant (where $\Gamma_o$ is the stabilizer of
$o$ in the automorphisms group $\aut(X)$ of the graph $X$) then
it is easy to prove that the generating functions $G$ and $F$
are $\Gamma_o$-invariant. Hence for any $\gamma \in \Gamma_o$,
the relation
$$
y \mapsto \lim_{n\rightarrow +\infty} \gamma(x_n)
$$
(where $x_n \rightarrow y \in \widehat X$) is a well defined map
which extends $\gamma$; this means that we can extend the subgroup
$\Gamma_o$ to a closed subgroup 
$\widetilde \Gamma_o$ of homeomorphic maps from $\widehat X$ onto itself. 
\hfill\break \noindent
It is quite obvious that, if $\mu$ is a finite variation measure on $X$,
${\cal F}$ a $\Gamma_o$-closed set of $|\mu|$-integrable maps (that is 
$f\gamma \in {\cal F}$ for any $f\in {\cal F}$ and for any $\gamma \in \Gamma_o$),
then any $f \in \Gamma$ is 
$\mu\circ\gamma$-integrable 
(for all $\gamma \in \Gamma$).
As a consequence, if $o\in X$ then for any 
$\gamma \in \Gamma_o$, TFAE
\begin{description}
\item{(i)} $\mu$ has the MVP with respect to ${\cal F}$ and $o$,
\item{(ii)} $\mu \circ \gamma$ has the MVP with respect to ${\cal F}$ and
$o$;
\end{description}
\noindent
in particular, the weak MVP property with respect to $o$ can be defined
for the equivalence classes of measure with respect to the equivalence 
relation $\mu \equiv_o \nu$ ($\mu \equiv_o \nu$ if and only if there exists $\gamma 
\in \Gamma_o$ such that $\mu=\nu\gamma$).
\hfill\break\noindent
Moreover the subset of all the regular boundary points is 
$\widetilde \Gamma_o$-invariant.
\end{REM}

\end{section}

\begin{section}[3]{Representation theorems for harmonic functions}
\label{3}

The topics in this section are essentially those of \cite{Cartier} with
the addition of some specific remarks and a result of 
\cite{Woess2}. This is indended to be a 
brief discussion about the main results of \cite{Cartier} which
we report here for sake of completeness and to fix the notations.

Let $(S,E(S))$ be a connected graph with an adapted (nearest neighbour) 
random walk $P$; chosen a reference vertex $o\in S$,
we denote by $\widehat S$ the Martin compactification of the random walk and
by $\partial S:=\widehat S \setminus S$ the space of ends. Moreover
let ${\cal F}$ the linear space of all the real-valued functions on $S$,
${\cal K}$ the linear subspace of ${\cal F}$ containing the functions
with finite support, ${\cal H}$ the linear subspace of ${\cal F}$ 
containing the harmonic functions and ${\cal D}$ the algebra of
the locally constant real-valued functions on $\widehat S$.

It is obvious that a function is locally constant if and only if
the counterimage of every set is open; in particular a locally constant
function is continuous and, if the domain is a compact topological space,
then its range is finite (this is our case). Note that the counterimage
$f^{-1}(A)$ of any set $A$ by means of any locally constant function
$f$ is open and closed.

If $(S,E(S))$ is a tree, $P$ an adapted transient random walk on $S$ and we
consider the Martin Kernels $k_o:S \times \widehat S \rightarrow \reali$,²
defined 
in the previous section, then the relation
$$
(K^*_o v)(x):= \sum_{y \in S} v(y)k_o(y,x), \qquad \forall x\in S, \forall v \in {\cal K}
$$
defines a linear map which turns out 
to be an isomorphism from ${\cal K}$ onto ${\cal D}$ (\cite{Cartier}, Proposition~A.1).
Differently stated, the set $\{k_{o,x}:x\in S\}$ (where for every $x \in S$ the function
$k_{o,x}$ is defined by $k_{o,x}(y):=k_o(x,y)$, $y\in \widehat S$) is a basis
of the linear space ${\cal D}$.

We call {\it distribution} on $\widehat S$ any linear functional on ${\cal D}$ and we
denote by ${\cal D}^*$ the set of distributions. We denote also by $\beta$
the set of all subsets of $\widehat S$ which are open and closed (this is trivially
a $\sigma$-algebra). The set $\{\chi_A:A \in \beta\}$ (where $\chi_A$ is the usual
characteristic function of a set $A$) is a set of generators
of ${\cal D}$; for any $\lambda \in {\cal D}^*$, we define a finitely additive
measure $\mu_\lambda$ on $\beta$ by $\mu(A):=\lambda(\chi_A)$ for every $A\in \beta$.
The correspondence $\lambda \mapsto \mu_\lambda$ is seen to be an isomorphism
from the space of distributions ${\cal D}^*$ onto the linear space 
of finitely additive (signed) measure on $\beta$.
This allows us to identify the distributions with the finite additive
measures, the action of $\mu_\lambda$ on $v \in {\cal D}$ is given by
the equality
$$
\lambda(v)=\int_{\widehat S} v \d \mu_\lambda = \sum_{\alpha \in {\rm Rg}(v)} \alpha 
\mu_\lambda(f^{-1}(\alpha)).
$$
We say that a distribution $\lambda$ is {\it supported} on a subset $E\subset \widehat S$
if and only if $\lambda(f)=0$ for any $f \in {\cal D}$ such that $f_{|E} \equiv 0$.

Moreover the following relations hold for any distribution $\lambda$ on 
$\widehat S$ and for any closed subset $E \subset \widehat S$,
$$
\begin{array}{rl}
\lambda(f) \geq 0, \ \forall f\in {\cal D}
\hbox{ such that } f\geq 0 
\Longleftrightarrow 
\mu_\lambda(A) \geq 0, \ \forall A\in \beta & \\
f \hbox{ is supported on } E 
\Longleftrightarrow 
\mu_\lambda(A) \geq 0, \ \forall A\in \beta, A\cap E =\emptyset; &
\end{array}
$$
in the first case we write $\lambda \geq 0$.

\begin{REM}\label{support}
If $E\subset \widehat S$, then $f_{|E} \equiv 0$ if and only if
$\lambda(f)=0$ for any distribution $\lambda$ supported on $E$. \hfill\break\noindent
Indeed the ``only if''
part is by definition; on the other hand, the Dirac measure $\delta_\xi$ is a $\sigma$-additive measure on ${\cal P}(\widehat S)$,
the restriction to $\beta$ can be identified with the distribution
$\lambda_\xi$ defined by $\lambda_\xi(f):=f(\xi)$ for any $\xi \in \widehat S$. 
Hence for any $\xi \in E$, $\lambda_\xi$ is supported on $E$ and
$0=\lambda_\xi(f)=f(\xi)$.
\end{REM}

We conclude this section stating, without any proof, 
the Poisson Representation Theorem for harmonic
functions on trees (\cite{Cartier}, Proposition~A.4) and discussing
the Poisson representation for bounded harmonic functions on general
set with general irreducible transition probabilities (\cite{Woess2},
Theorem 24.12).

\begin{thm}\label{representation}
Let $(S,P)$ be an irreducible, 
transient random walk on a tree; if $o \in S$
and $g \in {\cal F}$ then
\begin{description}
\item{(i)} there exists a unique distribution $\lambda_g$ on $\widehat S$ such that
$$
g(x)=\lambda(k_o(x,\cdot)), \ \forall x \in S;
$$
\item{(ii)} $g \in {\cal H}$ if and only if $\lambda_g$ is supported on 
$\partial X$;
\item{(iii)} $g$ is superharmonic and positive if and only if 
$\lambda_g \geq 0$.
\end{description}
\end{thm}

We call the unique distribution $\lambda_g$ the {\it Martin potential} of $g$.
An explicit expression for the map $I:{\cal H} \rightarrow {\cal D}^*$ such that
$I(g)=\lambda_g$ was found in \cite{Pica-Taible} (Theorem 2).

\begin{thm}\label{representation2}
Let $(X,P)$ be
a general irreducible, transient random
walk then
the equation
\begin{equation}\label{raphar2}
h(x):=\int_{\cal M}\varphi(\xi) k_o(x, \xi) \d \nu_o(\xi)
\equiv
\int_{\cal M} \varphi(\xi) \d \nu_x(\xi), \qquad \forall x \in X
\end{equation}
defines a harmonic function;
in particular 
the map $\varphi \mapsto h$, given by equation~\ref{raphar2}, is 
a linear, bicontinuous
operator from $L^\infty({\cal M}, \nu_o)$ onto ${\cal H}^\infty(X,P)$
\end{thm}


\end{section}

\begin{section}[5]{The case of measures with finite support on general graphs}
\label{5}

In this paragraph we characterize the weak MVP by means of the MVP
with respect to the set of harmonic functions
$\{k_o(\cdot,\xi)\}_{\xi \in \partial X}$.

We know from \cite{Knapp1} that the set of regular point
of the boundary is a Borel set, hence it is measurable. Also in the
transient case it could be an empty set; the Dirichlet problem is
solvable if and only if every boundary point is regular
(for the definition of
solvability of the Dirichlet problem see \cite{Woess2}, Chapter 4 
Paragraph 20). In the next theorem we are mostly interested in 
the representation of bounded harmonic functions by means of
the Martin compactification, for this reason it is crucial
that the set of regular points is ``sufficiently big'' (not
neccessarily the whole boundary).
When we say that a property is a.e.~true on $\partial X$ we mean w.~r.~to some
harmonic measure (that is the same, w.~r.~to every harmonic measure).

\begin{thm}\label{equiv2}
Let $(X,P)$ be an irreducible, transient random
walk, $o \in X$ and $\mu$ a (signed) measure with finite support.
Consider the following assertions
\begin{description}
\item{(i)} $\mu$ has the weak MVP with respect to $o$;
\item{(ii)} $L(k_o(\cdot,\xi),\mu)(o)=0$ a.~e.~on 
$\partial X$;
\end{description}
\noindent
then (ii) $\Longrightarrow$ (i). In particular
if the Martin compactification is almost regular
then (i) $\Longleftrightarrow$ (ii).
\end{thm}

\begin{pf}
(ii) $\Longrightarrow$ (i). It is Corollary~3.6 of \cite{Zucca1}.

(i) $\Longrightarrow$ (ii). If $\xi_0$ is a regular boundary point
then, defining
\begin{equation}\label{deffn}
\begin{array}{rl}
f_n(x):= \int_{\partial X} k_o(x_n, \xi)k_o(x,\xi)  \d 
\nu_o(\xi) \equiv &
\int_{\partial X} k_o(x,\xi) \d 
\nu_{x_n}, \\
& \qquad \forall x\in X, \forall n \in \naturali,
\end{array}
\end{equation}

where $\{x_n\}$ is a sequence in $X$ converging
to $\xi_0$ in the topology of $\widehat X$, we have that
$$
f_n(x) \ {\buildrel n \rightarrow +\infty \over \longrightarrow} \ 
k(x,\xi_0), \quad \forall x\in X
$$
because of Theorem~20.3 of \cite{Woess2}.
Now $f_n \in {\cal H}^\infty(X,P)$ for any $n\in \naturali$ and the support of $\mu$ is finite,
thus we have
$$
L(k_o(\cdot,\xi_0),\mu)(o)=\lim_{n\rightarrow +\infty}
L_o(f_n,\mu)(o) =0.
$$
\QED
\end{pf}

The previous theorem applies, in particular, if
$(X,P)$ is a nearest neighbour, transient random walk on a tree $X$ such that
the Martin boundary contains at least two elements: in that case
the set of regular point has $\nu_o$-measure 1 (see \cite{Woess4}, Theorem 4.2).
The previous theorem extends 
the equivalence between (i) and (iii) of Proposition 4.13 of \cite{Zucca1}. 

\end{section}

\begin{section}[4]{The case of measures with finite support on trees}
\label{4}

In this section we use of the Representation 
Theorem~\ref{representation} we stated in Section~\ref{3}. Our first aim is
to extend the results of \cite{Zucca1}, regarding the
finite variation measures with finite support on 
radial trees, to more general
settings (i.e.~general irreducible random walks on non-oriented trees).

\begin{thm}\label{equiv1}
Let $(T,P)$ an irreducible, transient random walk on a tree
with root $o$ and let $\mu$ a (signed) measure on $T$ with
finite support. Consider the following assertions:
\begin{description}
\item{(i)} $\mu$ has the weak MVP with respect to $o$,
\item{(ii)} $\mu$ has the strong MVP with respect to $o$,
\item{(iii)} for any $\xi \in \partial T$ we have that
$L(k_o(\cdot, \xi),\mu)(o)=0$;
\end{description}
\noindent
then (i) $\Leftarrow$ (ii) $\Longleftrightarrow$ (iii).
\end{thm}

\begin{pf}
(ii) $\Longleftrightarrow$ (iii). We easily note that from 
Theorem~\ref{representation}, (ii) is equivalent to
$$
\sum_{x \in \supp(\mu)} <\lambda, k_o(x,\cdot)>\mu(x) = 
\sum_{x \in \supp(\mu)} <\lambda, k_o(o,\cdot)>\mu(x)
$$
for any distribution $\lambda$ on $\widehat T$ with support in
$\partial T$. Using the finiteness of $\supp(\mu)$
and since the space of distribution separates the points
(according to Remark~\ref{support}), then
the previous relation is equivalent to 
$$
\sum_{x \in T} \mu(x) k_o(x,\cdot) = \mu(T) k_o(o,\cdot)
$$
that is 
$$
L(k_o(\cdot,\xi),\mu)(o) = 0, \quad \forall \xi \in \widehat T
\setminus T.
$$

(ii) $\Longrightarrow$ (i). It is trivial.
\end{pf}
\QED

We proved it in \cite{Zucca1} (Proposition~4.13)
that, under the condition of radiality, (i),(ii) and (iii) are equivalent.
Before extending this result to a bigger class of trees (see Corollary~\ref{trees1}) we give some examples
which prove that the equivalence between (i) and (ii) does not hold in general
(even if $\mu$ has finite support and $T$ is a tree or a transient
graph). 

\begin{exmp}\label{ex1}
We construct an example of a finite supported measure on a transient tree which
has the weak mean value property but not the strong one. Let us
take any tree $T_1$ (which has transient simple random walk) and a root $o$,
and let us consider now the tree $T$ obtained attaching to $T_1$ a straigh line 
(a copy of $\naturali$)
to the root $o$ by indentifying $o$ with the root $0 \in \naturali$
and let us call $P$ the simple random walk on $T$. 
It is simple to note that given any pair of distinct verticex $x,y$ on the straight
line (different from $o$) then for any $f\in {\cal H}(T,P)$ we have
$f_{|\naturali} \equiv const$ if and only if $f(x)=f(y)$. It is also 
easy to show that any bounded harmonic function $f$ must be constant 
on $\naturali$.
On such a tree is not difficult to construct a harmonic function
which is not constant (and hence not bounded) on $\naturali$: let $f_0$ one of those 
functions.
Let now consider the measure $\mu$ defined by 
$$
\mu(x):=\cases{
1 & if $x=1$\cr
-1 & if $x=2$ \cr
0 & if $x \not \in \{1,2\}$;\cr
}
$$
such a measure has the weak MVP with respect to any point
(since $\mu(T)=0$ and for any bounded harmonic function $f$ we have
$f(1)=f(2)$), but it cannot have
the strong MVP with respect to any vertex since 
$L(f_0,\mu)(x)= f_0(1)-f_0(2)\not = 0$.
Note that $T$ contains a recurrent branching subree, namely the copy of $\naturali$ (compare
with Corollary~\ref{trees1}).
\end{exmp}

\begin{exmp}\label{ex2}
Another example is provided by the simple random walk on the $d$-dimensional lattice 
$\zzz^d$: in this case the {weak Liouville property} holds
(every bounded harmonic function is constant), meanwhile
$\{f: \exists \alpha, \beta \in \reali: 
f(i_1, \ldots, i_d) \equiv \alpha i_1 +\beta\} \subset {\cal H}(\zzz^d,P)$.
Whence if $f \in {\cal H}(\zzz^d,P)$ is a nonconstant fuction we
cannot find a sequence $\{f_n\}$ in ${\cal H}^\infty(\zzz^d,P)$ converging
to $f$.
\end{exmp}

The equivalence between weak and strong MVP can be proved under certain conditions.

\begin{cor}\label{trees1}
Let $(T,P)$ be an irreducible, transient random walk on a tree, $o\in T$. 
Then TFAE:
\begin{description}
\item{(i)} every infinite branching subtree is transient;
\item{(ii)} for any $x,y \in T$ such that $x \sim y$ and $T_{x,y}$ is infinite we have $F(y,x)<1$;
\item{(iii)} there exists $x\in T$ ($\Leftrightarrow \forall x \in T$) such that
$\supp(\nu_x)=\partial T$;
\item{(iv)} there exists $x\in T$ ($\Leftrightarrow \forall x \in T$) such that
there exists a flux from $x$ to infinity, supported on each edge $(z,y)$ such that $T_{z,y}$ is infinite.
\end{description}
\noindent If the Martin boundary of $T$ contains at least two elements, then
any of the previous conditions is equivalent the following one:
\begin{description}
\item{(v)} given any signed measure with finite support $\mu$, then
$\mu$ has the weak MVP with respect 
to $o$ if and only if it has the strong MVP with respect to $o$
\end{description}
\noindent
\end{cor}
\begin{pf}
The equivalence between $(i)$, $(ii)$ and $(iii)$ follows by Theorem 7.5 of \cite{Woess3}.
The equivalence between one of the previous with $(iv)$ is straighforward
(remember that  any random walk on a tree is reversible).

(iii) $\Longrightarrow$ (v). 
Let $\supp(\nu_o)=\partial T$. If $\xi \in \partial T$ and $x\in B(o,n)$
such that 
$\supp(\mu)\subseteq B(o,n)$ and $x \wedge_o\xi =x$ then $\nu_o(\partial T_{o,x}) >0$.
By Theorem~\ref{equiv2} there exists $\xi_1\in \partial T_{o,x}$ such that
$L(k_o(\cdot,\xi_1),\mu)(o)=0$.
Since $k(y,\xi_1)=k(y,\xi_1 \wedge_o y)=k(y,\xi \wedge_o y)=k(y,\xi)$ for every 
$y \in T\setminus T_{o,x} \supseteq \supp(\mu)$ then we have
$k(\cdot,\xi_1)=k(\cdot,\xi)$ $\mu$ $a.~.e.~$ and hence
$$
0=L(k_o(\cdot,\xi_1),\mu)(o)=L(k_o(\cdot,\xi),\mu)(o).
$$
By Theorem~\ref{equiv1}, $(v)$ is proved. \hfill\break\noindent

(v) $\longrightarrow$ (i).
If $T_0$ is a recurrent, infinite branching subtree of 
the transient tree $T$, hence $\nu_x(\partial T_0)=0$ and every 
bounded harmonic function $f$ is constant on $T_0$; now
the measure $\nu$ in the last statement can be constructed as in Example~\ref{ex1}.
\end{pf}
\QED

The previous corollary applies, in particular, if 
$T$ is a locally finite tree
with minimum degree 2 and with 
finite upper bound to the lenghts of its unbranched geodesics 
(and $P$ is the simple random
walk): in this case the same property holds for any 
branching subree (and this implies easily
$(i)$). 

Moreover it is an extension of Proposition 4.13 of \cite{Zucca1} since for any irreducible, 
$\Gamma_o$-invariant random walk (see \cite{Zucca1} Definition 2.2) on an $o$-radial tree
the property (iii) of the previous corollary obviously holds.

\end{section}

\begin{section}[6]{The case of measures with unbounded support on general graphs}
\label{6}

If we want to deal with measures with countable support then we must slightly 
modify the hypotheses of Theorem~\ref{equiv2}.

\begin{thm}\label{equiv3}
Let $(X,P)$ be an irreducible, transient random
walk and $\mu$ a (signed) measure.
\begin{description}
\item{a)} For any $x,y \in X$ we have that $\sup_{z\in \widehat X}k_x(y,z)=\max_{z\in X}
k_x(y,z) = 1/F(x,y)$ and for any $x,y\in X$,
\begin{equation}\label{equiv3.1}
\sup_{z \in X} k_x(\cdot,z) \in L^1(|\mu|) \ \Longleftrightarrow 
\sup_{z \in X} k_y(\cdot,z) \in L^1(|\mu|). 
\end{equation}
\item{b)} Let $o\in X$ and $\sup_{z \in X} k_o(\cdot,z) \in L^1(|\mu|)$ (i.e.
$1/F(o,\cdot) \in L^1(|\mu|)$) then TFAE
\begin{description}
\item{(i)} $\mu$ has the weak MVP with respect to $o$;
\item{(ii)} $L(k_o(\cdot,\xi),\mu)(o)=0$ a.~e.~on 
$\partial X$;
\end{description}
then (ii) $\Longrightarrow$ (i). In particular
if the compactification is almost regular
then (i) $\Longleftrightarrow$ (ii).
\end{description}
\end{thm}

\begin{pf}
(a) By definition $k_x(y,z)=F(y,z)/F(x,z)$, since $F(x,z) \geq F(x,y)F(y,z)$ and
$F(x,y) \leq 1$, $F(x,x)=1$ then by continuity of $k_x$ we have that
\begin{equation}\label{equiv3.2}
\sup_{z\in \widehat X}k_x(y,z)=\max_{z\in X} k_x(y,z) = 1/F(x,y).
\end{equation}
Moreover, using again the previous inequality for $F$, we have that, for
any $x,y\in X$,
$$ 
{1 \over F(x,z)} \leq {1\over F(y,z)}{1 \over F(x,y)}, \qquad \forall z\in X
$$
and hence, by equation~\ref{equiv3.2}, we obtain equation~\ref{equiv3.1}.
\hfill\break\noindent
(b) \hfill\break
(ii) $\Longrightarrow$ (i). It is again Corollary~3.6 of \cite{Zucca1}.
\hfill\break
(i) $\Longrightarrow$ (ii). The proof is quite similar to the analogous case
of Theorem~\ref{equiv2} and we just outline the differences.
Since $\sup_{y\in X} k_o(x,y) \equiv
\sup_{\xi \in \widehat X} k_o(x,\xi) <+\infty$
by the Lebesgue Bounded Theorem (or just using uniform
convergence)
we have that
$\xi \mapsto \int_{X} k_o(\cdot,\xi) \d \mu$ is continuous on $\widehat X$
(and hence bounded by compactness) and the same holds for 
$\xi \mapsto \int_{X} k_o(\cdot,\xi) \d |\mu|$.
If $\xi_0$ is a regular boundary point $f_n$ is the sequence of functions defined by equation~\ref{deffn} from a sequence $\{x_n\}$ converging to $\xi_0$,
applying Fubini's Theorem we have that 
$k_o \in L^1(\nu_x \times |\mu|)$ and hence, using the Jordan decomposition of
$\mu$, we obtain
$$
L(f_n,\mu)(o)= \int_{{\cal M}} L(k_o(\cdot,\xi),\mu)(o) \d \nu_{x_n}.
$$
Now by the previous arguments $L(k_o(\cdot,\xi),\mu)(o)$ is continuous on 
$\partial X$ hence Theorem 20.3 of \cite{Woess2} yields again the conclusion.
\QED
\end{pf}

It is important to note the difference between (iii) of Theorem~\ref{equiv1}
and (ii) of Theorems~\ref{equiv2} and \ref{equiv3}; in
particular if $T$ is as in Example \ref{ex1}, then 
the end point of the straight line $\naturali$ is not a 
regular point.


\end{section}

\begin{section}[7]{Some results for general compactifications}
\label{7}

Throughout this whole section $\widehat X$ will always be
a natural compactification and $\{\nu_x\}$ the relative family of
harmonic measures. 
We state and prove two results which generalize those of the previous sections. In the first
Theorem we give a sufficient condition for the MVP with respect a particular family of functions:
in the case of the Martin compactification this is just the weak MVP.

\begin{thm}\label{7.1}
Let
$o\in X$ and ${\cal F}:=\{f:X \rightarrow \reali\ : f(x)=
\int_{\partial X} \phi \d \nu_x, \forall x \in X, \forall \phi \in L^\infty(\partial X,\nu_o)\}
$ and $\mu$ a finite variation (signed) measure on $X$ such that
$L(\nu_\cdot(A),\mu)(o)=0$ for every borel set $A \subset \partial X$, then
$\mu$ has the MVP with respect to ${\cal F}$ and $o$.
\end{thm}

\begin{pf}
We firstly note that, if $\phi \in L^\infty(\partial X, \nu_o)$, then 
$h_\phi(x):=\int_{\partial x} f \d \nu_x$ satisfies $h_\phi \in l^\infty(X)$ and
hence $h_\phi \in L^1(|\mu|)$.
If $s$ is a measurable simple function, $s=\sum_{i=1}^n \alpha_i \chi_{A_i}$, then
$$
L(h_s,\mu)(o)= L \left(\sum_{i=1}^\infty \alpha_i 
\nu_\cdot(A_i),\mu \right)(o)=\sum_{i=1}^\infty \alpha_i L(\nu_\cdot
(A_i),\mu)(o)=0.
$$
Let $\phi \in L^\infty(\partial X,\nu_o)$ and 
$\{s_n\}_n$ a sequence of measurable simple functions such that $s_n \rightarrow \phi$
and $|s_n| \leq |\phi|$, for any $n\in \naturali$, 
hence using the Bounded Convergence Theorem (on
the measurable space $(\partial X, \nu_x)$), we have
$$
h_{s_n}(x) :=\int_{\partial X} s_n \d \nu_x 
\rightarrow \int_{\partial X} \phi \d \nu_x=:h_\phi(x), \quad \forall x\in X.
$$
Moreover since $|h_{s_n}| \leq h_{|\phi|}$, according to Bounded
Convergence
Theorem again (applied to $(X,\mu)$), we have
$$
L(h_\phi,\mu)(o)=\lim_{n \rightarrow +\infty} L(h_{s_n},\mu)(o) =0.
$$
\QED
\end{pf}

In order to be able to guarantee the mean property of $\nu_\cdot(A)$ of the
previous Theorem one can check for the mean value property with respect
to a subfamily of ${\cal F}$.

\begin{thm}\label{7.2}
Let
$o\in X$, define
${\cal F}_1:=\{f:X \rightarrow \reali\ : f(x)=
\int_{\partial X} \phi \d \nu_x, \forall x \in X, \forall \phi \in C(\partial X)\}
$ and consider
$\mu$ a finite variation (signed) measure on $X$ which has the MVP
with respect to ${\cal F}_1$ and $o$ then
$L(\nu_\cdot(A),\mu)(o)=0$ for every borel set $A \subset \partial X$.
\end{thm}

\begin{pf}
Since $\partial X$ is a Hausdorff, compact topological space
and $\nu_o$ is regular
(according to Theorem~2.18 of \cite{Rudin1}),
we have that $C_c(\partial X)\equiv C(\partial X)$
is a dense subset of $L^p(\partial X, \nu_o)$ (see \cite{Rudin1} Theorem~3.14).
Let us choose a function $g \in l^1(X)$, $g(x)>0$ for any $x\in X$, and consider the measure
$\overline \nu$ defined by $\overline \nu(A):=\sum_{x \in X} g(x) \nu_x(A)$ for every
$A \subset \partial X$ borel subset.
If $\mu$ satisfies the MVP, from the density property for any 
borel set $A \subset
\partial X$,
there exists a sequence $\{\phi_n\}_n$ of continuous functions converging to 
$\chi_A$ in $L^1(\partial X,\overline \nu)$ (and hence in $L^1(\partial X, \nu_x)$
for any $x\in X$)
with the property $0\leq\phi_n\leq 1$. 
Moreover by Bounded Convergence Theorem (applied to $(X,\mu^+)$ and $(X,\mu^-)$)
since $\int_{\partial X} \phi_n\d \nu_x \rightarrow \nu_x(A)$, bounded by
$1$ (which is $|\mu|$-integrable since it is a finite variation measure), then
$$
L(\nu_\cdot(A), \mu)(o)=\lim_{n \rightarrow +\infty} L(\phi_n, \mu)(o)=0.
$$
\QED
\end{pf}

If ${\cal F}$ and ${\cal F}_1$ are as in Theorems~\ref{7.1} and \ref{7.2} 
respectively and $o\in X$ then the MVP with respect to ${\cal F}_1$ and $o$ implies
the MVP with respect to ${\cal F}$ and $o$.
This means that in the case of the Martin compactification the weak MVP
(and also the strong MVP if the hypotheses of Corollary~\ref{trees1} are
satisfied) is equivalent to the MVP with respect to the family ${\cal F}_1$.


Moreover if $(X,P)$ is a transient tree then the arguments of Section~\ref{2} imply
that the hypothesys of integrability of Theorem~{7.2} is equivalent to
$1/F(o, \cdot) \in L^1(|\mu|)$ and this is true for some $o\in X$ if and only
if it is true for any $o\in X$.

\end{section}
\end{article}

\end{document}